\documentclass[12pt]{article}
\usepackage[cp850]{inputenc}
\usepackage[T1]{fontenc}\usepackage{amsfonts}
\usepackage{graphicx}
\usepackage{amsmath}
\def\<{\langle}
\def\>{\rangle}
\date{}

\title{   Perpetuity property of the Dirichlet distribution}
\author{Pawe\l \ Hitczenko\thanks{Department of Mathematics, Drexel University, Philadelphia P.A. 19104,
 U.S.A. Partially supported by a grant from Simons Foundation (grant \#208766 to Pawe{\l} Hitczenko).}, G\'erard Letac
\thanks{Laboratoire de
Statistique et Probabilit\'es, Universit\'e Paul Sabatier, 31062 Toulouse,
France. This author thanks the Fields Institute for its hospitality during the preparation of this paper.}}
\begin{document}\maketitle
\begin{abstract} Let $X$,  $B$ and $Y$  be three Dirichlet,  Bernoulli and beta independent random variables such that $X\sim \mathcal{D}(a_0,\ldots,a_d),$ such that  $\Pr(B=(0,\ldots,0,1,0,\ldots,0))=a_i/a$ with $a=\sum_{i=0}^da_i$ and such that $Y\sim \beta(1,a).$ We prove that 
$X\sim X(1-Y)+BY.$ This gives the stationary distribution of a simple Markov chain on a tetrahedron.  We also extend this result to the case when $B$ follows  a quasi Bernoulli distribution $\mathcal{B}_k(a_0,\ldots,a_d)$ on the tetrahedron and when  $Y\sim \beta(k,a)$. We extend it even more generally to the case  where $X$ is a Dirichlet process and $B$ is a quasi Bernoulli random probability. Finally the case where the integer $k$ is replaced by a positive number $c$ is considered when $a_0=\ldots=a_d=1.$


\end{abstract}

\section{Introduction} In a recent paper [1],  Ambrus, Kevei and V\'igh make the following interesting observation: If $V, Y, W$ are independent random variables such that $V\sim \frac{1}{\pi}(\frac{1}{4}-v^2)^{-1/2}\textbf{1}_{(-1/2,1/2)}(v)dv$,  $Y$ is uniform on $(0,1)$ and  $\Pr(W=1)=\Pr(W=-1)=1/2$ then $$V\sim V(1-Y)+\frac{W}{2}Y.$$
The law $\mu$ of a random variable $V$ satisfying $V\sim VM+Q$ where the pair $(M,Q)$ is independent of $V$ on the right-hand side is often called a perpetuity generated by the law $\nu$ of $(M,Q)$. Thus, another way of stating the observation from [1] is that an arcsine random variable on $(-1/2,1/2)$ is a perpetuity  generated by the distribution of  $(M,Q)\sim(1-Y,WY/2)$. Part of the reason we found it interesting is that there are relatively few examples of exact solutions to perpetuity equations in the literature.  We will generalize this result, and our generalization will provide more examples of explicit generation of perpetuities, including  power semicircle distributions (see [2]). Note that the same perpetuity can be generated by different $\nu.$ See Section 6 below for a reminder of the classical link between perpetuities and the stationary distributions of the Markov chains obtained by iteration on random affine maps.  

 To carry out the generalization, we reformulate on $(0,1)$ the result in [1] by writing $X=V+\frac{1}{2}$ and $B=(1+W)/2.$ Clearly $X,Y,B$ are independent with $X\sim \beta(1/2,1/2)$, $B\sim \frac{1}{2}(\delta_0+\delta_1)$ and 
$X\sim X(1-Y)+BY.$ In Theorem 1.1 below, we give a generalization of the result in [1] expressed in terms of $X,Y$ and $B$. We will need the following notation. The natural basis of $\mathbb{R}^{d+1}$  is  denoted by $e_0,\ldots e_d$.  The convex hull of $\{e_0,\ldots e_d\}$ is a tetrahedron that we denote by $E_{d+1}.$ The elements of $E_{d+1}$ are therefore the vectors $\lambda=(\lambda_0,\ldots,\lambda_d)$ of $\mathbb{R}^{d+1}$ such that $\lambda_i\geq 0$ for $i=0,\ldots,d$ and such that $\lambda_0+\cdots+\lambda_d=1.$ 
If $p_0,\ldots,p_d$ are positive numbers with sum equal to one, the  distribution $\sum_{i=0}^dp_i\delta_{e_i}$ of $B=(B_0,\ldots, B_d)\in E_{d+1}$  is called a Bernoulli distribution. By definition $B$ satisfies $\Pr(B=e_i)=p_i.$  If $a_0,\ldots,a_d$ are positive numbers the Dirichlet distribution $\mathcal{D}(a_0,\ldots,a_d)$ of $X=(X_0,\ldots,X_d)\in E_{d+1}$ is such that the law  of $(X_1,\ldots,X_d)$ is
$$\frac{1}{B(a_0,\ldots,a_d)}(1-x_1-\cdots-x_d)^{a_0-1}x_1^{a_1-1}\ldots x_d^{a_d-1}\textbf{1}_{T_d}(x_1,\ldots,x_d)dx_1,\ldots dx_d$$ where $B(a_0,\ldots,a_d)=\frac{\Gamma(a_0)\ldots\Gamma(a_d)}{\Gamma(a_0+\ldots+a_d)}$ and 
where $T_d$ is the set of $(x_1,\ldots,x_d)$ such that $x_i>0$ for all $i=0,1,\ldots,d$, with the convention $x_0=1-x_1\cdots-x_d.$ For instance, if the real random variable $X_1$ follows the beta distribution $$\beta(a_1,a_0)(dx)=\frac{1}{B(a_1,a_0)}x^{a_1-1}(1-x)^{a_0-1}\textbf{1}_{(0,1)}(x)dx$$ then $(1-X_1,X_1)\sim \mathcal{D}(a_0,a_1).$ 

\vspace{4mm}\noindent \textbf{Theorem 1.1:} Let $a_0,\ldots,a_d$ be positive numbers. Denote $a=a_0+\cdots+a_d.$  
 Let $X$,  $Y$ and $B$  be three Dirichlet, beta and Bernoulli independent random variables such that $X\sim \mathcal{D}(a_0,\ldots,a_d)$ and $B\sim \sum_{i=0}^d\frac{a_i}{a}\delta_{e_i}$ are valued in $\mathbb{R}^{d+1}$ and such that $Y\sim \beta(1,a).$ Then 
$X\sim X(1-Y)+BY.$ 

\vspace{4mm}\noindent \textbf{Comments.} Considering each coordinate, Theorem 1.1 says that for all $i=0,\ldots,d$ we have  $X_i\sim X_i(1-Y)+B_iY.$ Since $1=\sum_{i=0}^dX_i=\sum_{i=0}^dB_i$ the statement for $i=0$ is true if it is verified for $i=1,\ldots,d.$ For instance for $d=1$ a reformulation of Theorem 1.1 is

\vspace{4mm}
\noindent \textbf{Theorem 1.2:} Let $a_0,a_1>0.$ Let $X_1,Y,B_1$ be three independent random variables such that
$X_1\sim\beta(a_1,a_0),\ Y\sim \beta(1,a_0+a_1),\ B_1\sim\frac{a_0}{a_0+a_1}\delta_{0}+\frac{a_1}{a_0+a_1}\delta_1.$ Then $X_1\sim X_1(1-Y)+B_1Y.$

\vspace{4mm}
\noindent Therefore the initial remark contained in [1] is a particular case of Theorem 1.1 for $d=1$ and  $a_0=a_1=1/2.$ More generally, the case $d=1$ and $a_0=a_1$ covers the power semicircle distributions discussed in [2] (with $\theta=a_0-3/2$). In particular, $a_0=a_1=3/2$ is the classical semicircle distribution. Theorem 1.1 is proved in Section 2. 
On the other hand, we shall see Theorem 1.1 as a particular case of the more general Theorem 4.1 in Section 4 below. Stating Theorem 4.1 needs the introduction of a new distribution $\mathcal{B}_k(a_0,\ldots,a_d)$ on the tetrahedron $E_{d+1}.$ We call it a quasi Bernoulli distribution of order $k.$  It is concentrated on the faces of order less than $k$ in a way that we will made reasonably explicit in Section 3. We add here a  family of laws with interesting properties to the zoo of distributions on a tetrahedron.

Finally, one  can  prove Theorem 1.2 by showing $\mathbb{E}(X_1(1-Y)+B_1Y)^n)=\mathbb{E}(X_1^n)$ for all integers $n.$ Our proof of  Theorem  4.1  is somewhat linked 
to  this  method of moments.  For this proof, we need to introduce  the $T_{c}$ transform of a distribution on the tetrahedron $E_{d+1}$. This is done in Section 2. We will  prove several important properties of  the $T_{c}$ transform in Theorem 2.1. Theorem 5.1 extends Theorem 4.1 to random probability measures on an abstract space $\Omega$, where the Dirichlet distribution is replaced by the so called Dirichlet process governed by the positive measure $\alpha$ on $\Omega.$  Surprisingly, this construction of Section 5 uses the Ewens distribution. Section 6  gives a  standard application of the preceeding results to certain Markov chains on $E_{d+1}.$

\section {The $T_{c}$ transform of a distribution on the tetrahedron}
In the sequel if $f=(f_0,\ldots,f_d)$ and $x=(x_0,\ldots,x_d)$ are in $\mathbb{R}^{d+1}$ we write $\<f,x\>=\sum_{i=0}^df_ix_i$ and we denote $$U_{d+1}=\{f=(f_0,\ldots,f_d)\in   \mathbb{R}^{d+1}\ ;\ f_0>0,\ldots,f_d>0\}.$$
Let $X=(X_0,\ldots,X_d)$ be a random variable on $E_{d+1}.$  Let $c>0.$ The $T_c$ transform of $X$  is the following function on $U_{d+1}:$
$$T_c(X)(f)=\mathbb{E}(\<f,X\>^{-c}).$$
Its existence is clear from $T_c(X)(f)\leq (\min_i f_i)^{-c}<\infty$. It satisfies $T_c(X)(\lambda f)=\lambda^{-c}T_c(X)(f).$ The  explicit calculation of $T_c(X)$ is easy in some rare cases, including the Dirichlet case $\mathcal{D}(a_0,\ldots,a_d)$ when $c=a=a_0+\ldots+a_d$ and the Bernoulli case  $\sum_{i=0}^kp_i\delta_{e_i}$. In some sense, the present paper originated from an effort to compute  $T_c(X)$ when $X\sim \mathcal{D}(a_0,\ldots,a_d)$ and $c=a+k$ where $k$ is a positive integer. For $d=1$, knowing  the $T_c$ transform is equivalent to knowing the function $t\mapsto \mathbb{E} ((1-tX)^{-c})$ on $(-\infty,1)$ when $X$ is a random variable valued in $[0,1]$ since  $$T_c((1-X,X))(1,1-t)=\mathbb{E} ((1-tX)^{-c}).$$ The $T_c$ transform is a tool which is 
in general  better adapted to the study of distributions on the tetrahedron than the Laplace transform $\mathbb{E}(\exp(-\<f,X\>)).$ The next theorem gathers its main properties. It shows for instance that $T_c(X)$ characterizes the distribution of $X$ and gives in (\ref{FF1}) a crucial probabilistic interpretation to the product $T_a(X_0)T_b(X_1)$ when 
$X_0$ and $X_1$ are independent random variables valued in $E_{d+1}.$

\vspace{4mm}
\noindent \textbf{Theorem 2.1:} \begin{enumerate}
\item If $X$ and $Z$ are  random variables on $E_{d+1}$ and if there exists $c>0$ such that $T_c(X)(f)=T_c(Z)(f)$ for all $f\in U_{d+1}$ then $X\sim Z.$ 
\item If $k$ is a non--negative integer and if  $H=-(\frac{\partial }{\partial f_0}+\cdots+\frac{\partial }{\partial f_d})$ then 
\begin{equation}\label{DIFF}H^kT_c(X)=(c)_kT_{c+k}(X),\end{equation}
where $(c)_n$ is the Pochhammer symbol defined by $(c)_0=1$ and $(c)_{n+1}=(c)_n(c+n).$ 
\item If $(a_0,\ldots,a_d)\in U_{d+1}$ with $a=a_0+\ldots+a_d$ and if $X\sim \mathcal{D}(a_0,\ldots,a_d)$ then \begin{equation}\label{CLD}T_a(X)(f)=f_0^{-a_0}\ldots f_d^{-a_d}.\end{equation} 
\item Suppose that $X_0,\ldots, X_r, Y$ are independent random variables  such that $X_i\in E_{d+1}$ for $i=0,\ldots,r$ and $Y=(Y_0,\ldots,Y_r)\in E_{r+1}$ has Dirichlet distribution $\mathcal{D}(b_0,\ldots,b_r).$ Then for   $b=b_0+\ldots+b_r$ and for $Z=X_0Y_0+\ldots +X_rY_r$ we have on $ U_{d+1}:$
\begin{equation}\label{FF}
T_b(Z)(f)=T_{b_0}(X_0)(f)\ldots T_{b_d}(X_d)(f). 
\end{equation}
In particular, if $Y\sim \beta(b_1,b_0)$   we have \begin{equation}\label{FF1}T_{b_0+b_1} ((1-Y)X_0+YX_1)=T_{b_0}(X_0)T_{b_1}(X_1).\end{equation}
\item The probability of the face $x_0=\ldots=x_k=0$ is computable by the $T_c$ transform:
$$\lim_{f_0\rightarrow \infty}T_c(X)(f_0,\ldots,f_0,1,1,\ldots,1)=\Pr (X_0=X_1=\ldots=X_k=0).$$

\end{enumerate}

\vspace{4mm}
\noindent \textbf{Proof:} For part 1) fix $g\in \mathbb{R}^{d+1},$ set $f_i=1-tg_i$ for $t$ small enough and develop $t\mapsto \mathbb{E}(\<f,X\>^{-c})$ in a neighborhood of $t=0.$  Since $\<f,X\>=1-t\<g,X\>$ we have
$$T_c(X)(f)=\mathbb{E}((1-t\<g,X\>)^{-c})=\sum_{n=0}^{\infty}\frac{(c)_n}{n!}\mathbb{E}(\<g,X\>^n)t^n.$$
It follows from the hypothesis $T_c(X)=T_c(Z)$ that $\mathbb{E}(\<g,X\>^n)=\mathbb{E}(\<g,Z\>^n)$ for all $n.$ 
Thus $\<g,X\>\sim \<g,Z\>$ since both are bounded random variables with the same moments. Since this is true for all 
$g\in \mathbb{R}^{d+1}$ we have $X\sim Z.$ Formula (\ref{DIFF}) is easy to obtain by induction on $k$ using the fact that $X_0+\ldots+X_d=1.$ Let us give a proof of the standard formula (\ref{CLD}) by the so called beta-gamma algebra. It differs from the method of Proposition 2.1 in [4]. We write $\gamma_c(dv)=e^{-v}v^{c-1}\textbf{1}_{(0,\infty)}(v)dv/\Gamma(c).$ Consider independent $V_0,\ldots, V_d$ such that $V_i\sim \gamma_{a_i}$ and define $V=V_0+\ldots+V_d$ and $X_i=V_i/V$ for all $i=0,\ldots,d.$ Recall that $(X_0,\ldots,X_d)\sim \mathcal{D}(a_0,\ldots,a_d)$ is independent of $V\sim \gamma_a.$ Therefore
$$\mathbb{E}\left(\frac{1}{\<f,X\>^a}\right)=\mathbb{E}\left(\int_0^{\infty}e^{-v\<f,X\>}v^{a-1}\frac{dv}{\Gamma(a)}\right)=\mathbb{E}\left(e^{V-V\<f,X\>}\right)=\mathbb{E}\left(e^{\sum_{i=0}^d(V_i-f_iV_i}\right)=\prod_{i=0}^d\frac{1}{f_i^{a_i}}.$$
Formula (\ref{FF}) follows from  (\ref{CLD}) by replacing $X,a_0,\ldots,a_d$ by $Y,b_0,\ldots,b_r$ and $f$ by $\<f,X_0\>,\ldots,\<f,X_r\>.$  Using conditioning and the independence of $X_0,\ldots,X_r$ we obtain
$$T_b(Z)(f)=\mathbb{E}\left(\mathbb{E}([\sum_{j=0}^rY_j\<f,X_j\>]^{-b}|X_0,\ldots,X_r)\right)=\mathbb{E}\left(\prod_{j=0}^r\<f,X_j\>^{-b_j}\right)=\prod_{j=0}^r T_{b_j}(X_j)(f).$$
Applying (\ref{FF}) to $(Y_0,Y_1)=(1-Y,Y)\sim \mathcal{D}(b_0,b_1)$ we get $Z=(1-Y)X_0+YX_1.$ This leads to (\ref{FF1}). Property 5 is obvious since the events $X_0+\cdots+X_k=0$ and $X_0=X_1=\ldots=X_k=0$ coincide.  $\square$

\vspace{4mm}\noindent \textbf{Proof of Theorem 1.1} We prove Theorem 1.1 by taking $X_0=X,$ $X_1=B$, $b_1=1$ and $b_0=a$ in (\ref{FF1}). Thus 
$$T_1(B)(f)=\frac{1}{a}\left(\frac{a_0}{f_0}+\cdots+\frac{a_d}{f_d}\right).$$ The trick for computing $T_{1+a}(X)$ is to observe from (\ref{DIFF}) and (\ref{CLD}) that
\begin{eqnarray*}T_{1+a}(X)(f)&=&\frac{-1}{a}\left(\sum_{i=0}^d\frac{\partial}{\partial f_i}\right)\prod_{i=0}^d\frac{1}{f_i^{a_i}}=T_a(X)(f)T_1(B)(f).\end{eqnarray*}
From (\ref{FF1}) we also know that for $Z=(1-Y)X+YB$ we have  $T_{1+a}(Z)=T_a(X)T_1(B).$ Thus $T_{1+a}(Z)=T_{1+a}(X).$  Part 1 of Theorem 2.1 implies $X\sim Z.$ $\square$
\section{The quasi Bernoulli distributions on a tetrahedron}
First, we slightly extend the definition of a Dirichlet distribution $\mathcal{D}(a_0,\ldots,a_d)$ by allowing $a_i\geq 0$ instead of $a_i>0$ while keeping $a=\sum_{i=0}^da_i>0.$ 
For such a sequence $(a_0,\ldots,a_d)$ we define the nonempty set $T=\{i\ ; \ a_i>0\}$. We say that $\mathcal{D}(a_0,\ldots,a_d)$ is the Dirichlet distribution concentrated on the tetrahedron $E_T$ generated by $(e_i)_{i\in T}$ with parameters  $(a_i)_{i\in T}.$
If $X\sim \mathcal{D}(a_0,\ldots,a_d)$ the formula $\mathbb{E}(\<f,X\>^{-a})=\prod_{i=0}^df_i^{-a_i}$ still holds. If $T$ contains only one element $i_0$ then $\mathcal{D}(a_0,\ldots,a_d)$ is simply $\delta_{e_{i_0}}$ and does not depend on $a.$ 
Now recall a simple combinatorial formula where $k$ is a positive integer and $a=a_0+\cdots+a_d:$ 
\begin{equation}\label{WW}\sum_{(b_0,\ldots,b_d)\in \mathbb{N}^{d+1};\atop  b_0+\cdots+b_d=k}\prod_{i=0}^d\frac{(a_i)_{b_i}}{b_i!}=\frac{(a)_k}{k!}.\end{equation}
The proof is immediate if we use  generating functions: expand $\prod_{i=0}^d(1-t)^{-a_i}=(1-t)^{-a}$ in a power series on both sides. We now define our new distributions.

Let $a_0,\ldots,a_d>0$ and $a=a_0+\cdots+a_d$ and let $k$ be a positive integer. The quasi Bernoulli distribution of order $k$ is the distribution on the tetrahedron $E_{d+1}$ defined as the  mixing of Dirichlet distributions
\begin{equation}\label{QBk}
\mathcal{B}_k(a_0,\ldots,a_d)=\frac{k!}{(a)_k}\sum_{(b_0,\ldots,b_d)\in \mathbb{N}^{d+1};\atop  b_0+\cdots+b_d=k}\prod_{i=0}^d\frac{(a_i)_{b_i}}{b_i!}\mathcal{D}(b_0,\ldots,b_d).\end{equation}
Formula (\ref{WW}) shows that (\ref{QBk}) is indeed a probability on $E_{d+1}.$ Setting $c=k$ in Theorem 4.3 below gives an explicit form of $\mathcal{B}_k(a_0,\ldots,a_d)$ in the particular case $a_0=\ldots=a_d=1.$ 
For the sake of simplicity of the next statement denote \begin{equation}\label{SIG}\sigma_j=\sum_{i=0}^d\frac{a_i}{f_i^j}.\end{equation}

\vspace{4mm}
\noindent \textbf{Proposition 3.1:} If $B\sim \mathcal{B}_k(a_0,\ldots,a_d)$ then\begin{eqnarray}\label{QBE}
T_k(B)(f)&=&\frac{k!}{(a)_k}\sum_{(b_0,\ldots,b_d)\in \mathbb{N}^{d+1};\atop  b_0+\cdots+b_d=k}\prod_{i=0}^d\frac{(a_i)_{b_i}}{b_i!f_i^{b_i}}\\\label{QBEX}
&=&\frac{k!}{(a)_k}\sum_{(m_1,\ldots,m_k)\in \mathbb{N}^{k};\atop m_1+2m_2+\cdots+km_k=k}
\prod_{j=1}^k\frac{\sigma_j^{m_j}}{j^{m_j}m_j!}.
\end{eqnarray}

\vspace{4mm}
\noindent \textbf{Proof.} Formula (\ref{QBE}) is obvious from the definition of $\mathcal{B}_k(a_0,\ldots,a_d)$ and formula (\ref{CLD}). To prove (\ref{QBEX}) denote by $\frac{k!}{(a)_k}B_k$ and by $\frac{k!}{(a)_k}C_k$ the right--hand sides of (\ref{QBE}) and (\ref{QBEX}) respectively.  Now 
\begin{eqnarray*}\sum_{k=0}^{\infty}C_kt^k &=&\sum_{m_1,m_2,\ldots}\prod_{j\geq 1}\frac{t^{jm_j}\sigma_j^{m_j}}{j^{m_j}m_j!}=\exp\left(\sum_{j\geq 1}\frac{t^j\sigma_j}{j}\right)\\&=&\exp\left(\sum_{j\geq 1}\frac{t^j}{j}\sum_{i=0}^d\frac{a_i}{f_i^j}\right)=\prod_{i=0}^d\left(1-\frac{t}{f_i}\right)^{-a_i}.\end{eqnarray*}
Similarly one computes $\sum_{k=0}^{\infty}B_kt^k$, leading to $B_k=C_k$ and (\ref{QBEX}). $\square$
\vspace{4mm}
\noindent

The remainder of this section is made of several remarks on $\mathcal{B}_k(a_0,\ldots,a_d).$ Section 4 contains further information about it. If $T\subset\{0,\ldots,d\}$ denote by $F_T$ the relative interior of $E_T.$ This set is sometimes called a face of $E_{d+1}.$ It is equal to the relative interior of $E_{d+1}$ if $T=\{0,\ldots,d\}$ and the family of $F_T$'s is a partition of $E_{d+1}.$ Therefore 
$\mathcal{B}_k(a_0,\ldots,a_d)$ is a mixing of distributions on the faces $F_T$ which have densities $h_{k,T}$ with respect to the uniform distribution $\lambda_T$ on $F_T.$ Here $\lambda_T=\mathcal{D}(b_0,\ldots,b_d)$ where $b_i=1$ if $i\in T$ and $b_i=0$ if not. Observe that if $k\leq d$ only faces of dimension less than $k$ are charged by 
$\mathcal{B}_k(a_0,\ldots,a_d).$ To be more specific, denote $a_T=\sum_{i\in T}a_i$ and $b_T=\sum_{i\in T}b_i.$ 
When restricted to $(a_i)_{i\in T}$ formula (\ref{WW}) becomes 
$$\sum_{(b_i)_{i\in T}\in \mathbb{N}^{T};\  b_T=k}\prod_{i\in T}\frac{(a_i)_{b_i}}{b_i!}=\frac{(a_T)_k}{k!}.$$
A probabilistic interpretation of this is   $$\mathcal{B}_k(a_0,\ldots,a_d)\left(\bigcup_{S\subset T}F_S\right)=
(a_T)_k/(a)_k.$$ Since the $F_S$ are disjoint, 
for $S\subset \{0,\ldots,d\}$ the weights $w_S=\mathcal{B}_k(a_0,\ldots,a_d)(F_S)$  satisfy 
$\sum_{S\subset T}w_S=\frac{(a_T)_k}{(a)_k}.$ The principle of inclusion--exclusion therefore implies that
$w_T=\frac{1}{(a)_k}\sum_{S\subset T}(-1)^{|T\setminus S|}(a_S)_k.$ Let us introduce the symmetric  polynomial 
$$P_k(a_0,\ldots,a_d)=\sum_{S\subset\{0,\ldots,d\}}(-1)^{d+1-|S|}(a_S)_k.$$ Its explicit calculation is not easy. With the convention $P_0=1,$ here is a generating function:
\begin{eqnarray}\nonumber\sum_{k=0}^{\infty}P_k(a_0,\ldots,a_d)\frac{t^k}{k!}&=&\sum_{S\subset\{0,\ldots,d\}}(-1)^{d+1-|S|}\sum_{k=0}^{\infty}(a_S)_k
\frac{t^k}{k!}=\sum_{S\subset\{0,\ldots,d\}}(-1)^{d+1-|S|}(1-t)^{-a_i}\\&=&\label{FPk}\prod_{i=0}^d[(1-t)^{-a_i}-1]=t^{d+1}\prod_{i=0}^d\frac{(1-t)^{-a_i}-1}{t}.\end{eqnarray}
In particular, formula (\ref{FPk}) shows that $P_k(a_0,\ldots,a_d)=0$ if $k\leq d$ and that
$$P_{d+1}(a_0,\ldots,a_d)=(d+1)!a_0\ldots a_d.$$
With this notation we have $w_T=\frac{1}{(a)_k}P_k\left((a_i)_{i\in T}\right)$ (recall that $\sum_{T\subset\{0,\ldots,d\}}w_T=1).$ Another representation of the quasi Bernoulli distribution as a sum of mutually singular measures is  \begin{equation}\label{WQB}\mathcal{B}_k(a_0,\ldots,a_d) =\sum_{T\subset\{0,\ldots,d\}}w_T h_{k,T}\lambda_{T}.\end{equation} 
For simplicity denote $h_{k,T}$ by $h_{k,d}$ in the particular case $T=\{0,\ldots,d\}.$ Of course it is not zero only if $k\geq d+1.$ 
The following proposition gives a generating function for the sequence $(h_{k,d}(x))_{k\geq d+1}$ in terms of confluent functions

\vspace{4mm}
\noindent \textbf{Proposition 3.2:} For $a,b>0$ denote $_1F_1(a;b;z)=\sum_{n=0}^{\infty}\frac{(a)_n}{n!(b)_n}z^n.$ Then
$$\sum_{k=d+1}^{\infty}\frac{1}{(k-1)!}P_k(a_0,\ldots,a_d)h_{k,d}(x_0,\ldots,x_d)t^{k-d-1}=\prod_{i=0}^d\frac{1}{a_i}\ _1F_1(a_i+1\,;\,2\,;\,x_it).$$

\vspace{4mm}
\noindent \textbf{Proof:} Restricting (\ref{WQB}) to the interior of $E_{d+1}$ we get, by writing $n_i=b_i-1$ and by using the definition of the Dirichlet distribution
\begin{eqnarray*}P_k(a_0,\ldots,a_d) h_{k,d}(x)\textbf{1}_{E_{d+1}}(x)dx&=&\sum_{b_i>0\ \forall i; \atop \sum_{i=0}^db_i=k}\left(\prod_{i=0}^k\frac{(a_i)_{b_i}}{b_i!}\right)\mathcal{D}(b_0,\ldots,b_d)(dx)\\&=&(k-1)!
\sum_{n_i\geq 0\ \forall i;\ \atop \sum_{i=0}^dn_i=k-d-1}\left(\prod_{i=0}^d\frac{(a_i)_{n_i+1}}{(n_i+1)!n_i!}x_i^{n_i}\right)\textbf{1}_{E_{d+1}}(x)dx.\end{eqnarray*}
Multiplying both sides by $t^{k-d-1}$ and summing on $k=d+1,d+2,\ldots$ will give the proposition. $\square$

\section{Perpetuities for quasi Bernoulli} We compute now the $T_k$ transform of a quasi Bernoulli distribution $\mathcal{B}_k(a_0,\ldots,a_d)$ and we deduce from it the desired extension of Theorem 1.1.

\vspace{4mm}
\noindent \textbf{Theorem 4.1:} Let $a_0,\ldots,a_d>0$ with $a=a_0+\cdots+a_d$ and let $k$ be a positive integer. Suppose that  $X\sim \mathcal{D}(a_0,\ldots,a_d)$ and $B\sim\mathcal{B}_k(a_0,\ldots,a_d).$ Then 
$$T_k(B)(f)=\frac{T_{a+k}(X)(f)}{T_a(X)(f)}.$$ In particular, if $X,B$ and $Y\sim \beta(k,a)$ are independent then 
$$X\sim (1-Y)X+YB.$$

\vspace{4mm}
\noindent \textbf{Proof:} We have introduced the differential operator $H$ on $U_{d+1}$ in Theorem 2.2. 
 Consider the function $F(f)=\prod_{i=0}^d f_i^{-a_i}=T_a(X)(f).$ The idea of the proof is to compute $F^{-1}H^k (F)$ in two ways. A multinomial expansion shows that
$$H^{k}=k!\sum_{(b_0,\ldots,b_d)\in \mathbb{N}^{d+1};\atop  b_0+\cdots+b_d=k}\prod_{i=0}^d\frac{(-1)^{b_i}}{b_i!}\frac{\partial^{b_i}}{\partial f_i^{b_i}}.$$ We also observe that
$$\left(\prod_{i=0}^d(-1)^{b_i}\frac{\partial^{b_i}}{\partial f_i^{b_i}}\right)F=F\sum_{i=0}^d\frac{(a_i)_{b_i}}{ f_i^{a_i+b_i}}.$$ Combining these last two formulas with the definition of $\mathcal{B}_k(a_0,\ldots,a_d)$ we obtain that $F^{-1}H^k(F)=(a)_kT_k(B).$ On the other hand by applying formula (\ref{DIFF}) to $X\sim \mathcal{D}(a_0,\ldots,a_d)$
and to $c=a$ we get $F^{-1}H^k(F)=(a)_kT_{a+k}(X).$ Comparing the two results yields the proof of $T_a(X)T_k(B)(f)=T_{a+k}(X).$ Applying (\ref{FF1}) completes the proof of the theorem. $\square$

\vspace{4mm}
\noindent \textbf{Corollary 4.2 :} $\lim_{k\rightarrow\infty}\mathcal{B}_k(a_0,\ldots,a_d)=\mathcal{D}(a_0,\ldots,a_d)$ where the convergence is in the weak sense.

\vspace{4mm}
\noindent \textbf{Proof:} If $X\sim \mathcal{D}(a_0,\ldots,a_d)$, $Y_k\sim \beta(k,a)$ and $B_k\sim \mathcal{B}_k(a_0,\ldots,a_d)$ are independent, Theorem 4.1 shows that $(1-Y_k)X+Y_kB_k\sim \mathcal{D}(a_0,\ldots,a_d).$ Since $E_{d+1}$ is compact there exists a subsequence $k_n \rightarrow_{n\rightarrow\infty} 
\infty$ and a probability $\mu$ on $E_{d+1}$ such that $\mathcal{B}_{k_n}(a_0,\ldots,a_d)\rightarrow_{n\rightarrow\infty} \mu$ in the weak sense. Furthermore the distribution of $1-Y_k$ converges to the Dirac mass $\delta_0:$ a quick way to see this is to consider the Mellin transform for $s>0:$  $$\mathbb{E}((1-Y_k)^{s})=\frac{\Gamma(a+s)}{\Gamma(a)}\frac{\Gamma(a+k)}{\Gamma(a+k+s)}\rightarrow_{k\rightarrow\infty} 
0.$$Therefore the distribution of $(1-Y_k)X+Y_kB_k$ converges weakly to $\mu.$ As a consequence $\mu=\mathcal{D}(a_0,\ldots,a_d)$ and does not depend on the particular subsequence $(k_n).$ This proves the result. $\square$

\vspace{4mm}
\noindent  Theorem 4.1 implies that for all integer $k$ and for $X\sim \mathcal{D}(a_0,\ldots,a_d)$ there exists a probability distribution for $B$ such that $T_k(B)=\frac{T_{k+a}(X)}{T_a(X)}.$ One can wonder whether this statement can be extended to positive real numbers $c$ instead of the integers $k$ or not. More specifically, does there exist a distribution $\mathcal{B}_c(a_0,\ldots,a_d)$ on $E_{d+1}$ for $B$ such that $T_c(B)=\frac{T_{c+a}(X)}{T_a(X)}?$  We observe easily that it cannot be true: Taking  $c$ as a  positive number, and  $X$ uniform on $(0,1)$ and $Y\sim \beta(c,2)$ with $X$ and $Y$ independent; then \begin{itemize}\item If $0<c<1$ it is impossible to find a distribution for a random variable $B$ independent of $X,Y$ such that $X\sim (1-Y)X+YB.$\item If $c\geq 1$ and if $B\sim \frac{1}{c+1}(\delta_0+\delta_1)+\frac{c-1}{c+1}1_{(0,1)}(b)db$ is independent of $(X,Y)$ then $$X\sim (1-Y)X+YB.$$\end{itemize}
More generally we prove the following 

\vspace{4mm}
\noindent \textbf{Theorem 4.3.} Let $c$ be a  positive number. For a non--empty set $T\subset\{0,\ldots,d\}$ we denote by $\lambda_T$ the uniform probability on the convex set generated by 
$\{e_i\ ;\ i\in T\}.$ Introduce also the uniform probability on the union of the faces of $E_{d+1}$ of dimension $k$ $$\Lambda_k=\frac{(k+1)!(d-k)!}{(d+1)!}\sum_{T\subset \{0,\ldots,d\}; |T|=k+1}\lambda_T$$
and consider the signed measure on $E_{d+1}$ defined by
$$\nu_{c,d}=\frac{d!(d+1)!}{(c+1)(c+2)\ldots (c+d)}\sum_{k=0}^{d} \frac{(c-1)(c-2)\ldots (c-k)}{k!(k+1)!(d-k)!}\Lambda_k.$$
Then for all $f_i>0$, $i=0,\ldots,d:$ \begin{equation}\label{CP}\int_{E_{d+1}}\frac{\nu_{c,d}(dx)}{\<f,x\>^c}=f_0\ldots f_d\int_{E_{d+1}}\frac{\Lambda_d(dx)}{\<f,x\>^{c+d+1}}.\end{equation}
In particular,  if $Y\sim \beta(c,d+1)$ and $X\sim \Lambda_d$ are independent, 
then there exists a random variable $B$ on $E_{d+1}$ independent of $(Y,X)$ such that $X\sim (1-Y)X+YB$ if and only if either $c$ is a non--negative integer or if $c>d.$ Under these circumstances $B\sim \nu_{c,d}.$

\vspace{4mm}
\noindent \textbf{Comments.}  Note that $\Lambda_d=\mathcal{D}(1,\ldots,1).$ Therefore the theorem says that the quasi Bernoulli distribution $\nu_{c,d}=\mathcal{B}_c(1,\ldots,1)$ with continuous parameter $c$ does exist if and only if either $c$ is an integer or is $>d.$ For $d=2$   denote by $\lambda_{ij}$  the uniform distribution on the segment $e_i,e_j$, and by $\Lambda_2$ the uniform distribution of the triangle with vertices $e_0,e_1,e_2.$ Then for $c=1$ or $c\geq 2$ 
$$\nu_{c,2}=\mathcal{B}_c(1,1,1)=\frac{1}{(c+1)(c+2)}\left(2(\delta_{e_0}+\delta_{e_1}+\delta_{e_2})+2(c-1)(\lambda_{01}+\lambda_{02}+\lambda_{12})+(c-1)(c-2)\Lambda_2\right).$$
The proof of Theorem 4.3 is intricate enough to lead us to think that the existence of $\mathcal{B}_c(a_0,\ldots,a_d)$ for arbitrary positive numbers $(a_0,\ldots,a_d)$ is a delicate problem, even when all $a_i$ are equal. To illustrate this for $d=2$ we have to study whether there exists or not a positive probability $\mu(db)$ on $[0,1]$, depending on $a_0$, $a_1$ and $c,$ such that for all $(f_0,f_1)=(1,1-t)$ we have
$$\int_{0}^1\frac{\mu(db)}{(1-tb)^{c}}=\frac{(1-t)^{a_1}}{B(a_0,a_1)}\int_0^1\frac{x^{a_1-1}(1-x)^{a_0-1}}{(1-tx)^{a_0+a_1+c}}dx.$$ Application of Property 5 of Theorem 2.1 shows that necessarily $\mu(db)$ has atoms at $0$ and $1$. As shown by Proposition 4.4 below the mass at $0$ is given by the limit
\begin{equation}\label{IP}\lim_{t\rightarrow -\infty}\frac{(1-t)^{a_1}}{B(a_0,a_1)}\int_0^1\frac{x^{a_1-1}(1-x)^{a_0-1}}{(1-tx)^{a_0+a_1+c}}dx=\frac{B(a_0+c,a_1)}{B(a_0,a_1)}.\end{equation}
A similar result holds for the mass at 1.  However finding the part of $\mu(db)$ concentrated on $(0,1)$ is challenging. One can postulate that it has a density $f$ 
which therefore satisfies, in terms of the Gauss hypergeometric function $_2F_1:$ 
\begin{equation}\label{HE}\int_{0}^1\frac{f(b)db}{(1-tb)^c}=(1-t)^{a_1}\ _2F_1(a_0+a_1+c,a_1;a_0+a_1,t)-\frac{B(a_0+c,a_1)}{B(a_0,a_1)}-\frac{B(a_0,a_1+c)}{B(a_0,a_1)}\frac{1}{(1-t)^c}.\end{equation}
If $a_0$ and $a_1$ are positive integers one can show that $f$ is a polynomial of degree $a_0+a_1-2$ with a complicated expression.

\vspace{4mm}
\noindent \textbf{Proof of Theorem 4.3.}  Since all probabilities $\Lambda_0,\ldots,\Lambda_d$ are mutually singular, clearly $\nu_{c,d}$ is a positive measure if and only if either $c$ is an integer or if $c>d.$
We observe also that 
\begin{equation}\label{GF}\sum_{k=0}^d\frac{1}{k!(k+1)!(d-k)!}(c-1)(c-2)\ldots(c-k)=\frac{1}{d!(d+1)!}(c+1)(c+2)\ldots(c+k).\end{equation}
 (Compute the coefficient of $z^{d+1}$ on both sides of $(1+z)^{d}(1+z)^{c}=(1+z)^{d+c}$ to see this). This proves that the total mass of $\nu_{c,d}$ is one. 
 Denote for simplicity $F_c(f_0,\ldots,f_d)=\int_{E_{d+1}}\frac{\Lambda_d(dx)}{\<f,x\>^{c+d+1}}.$ This is a symmetric function of the $f_i's.$ We now show by induction on $d$ that \begin{equation}\label{BF}F_c(f_0,\ldots,f_d)=\frac{d!}{(c+1)(c+2)\ldots(c+d)}\sum_{i=0}^d\frac{1}{f_i^{c+1}\prod_{j\neq i}(f_j-f_i)}.\end{equation}This is correct for $d=1$ since 
$$\int_{0}^1\frac{dx_1}{(f_0(1-x_1)+f_1x_1)^{c+2}}=\frac{1}{(c+1)(f_0-f_1)}\left(\frac{1}{f_1^{c+1}}-\frac{1}{f_0^{c+1}}\right).$$Assuming that (\ref{BF}) is true for $d-1$ we write (recall that $T_d$ is the tetrahedron defined in Section 1 and that its Lebesgue measure is $1/d!$):
\begin{eqnarray*}&&F_c(f_0,\ldots,f_d)=d!\int_{T_d}\frac{dx_1\ldots dx_d}{(f_0(1-x_1-\cdots-x_{d})+f_1x_1+\cdots+f_dx_d)^{c+d+1}}\\&=& d!\int_{T_{d-1}}\left(\int_0^{1-x_1-\cdots-x_{d-1}}\frac{dx_d}{(f_0(1-x_1-\cdots-x_{d})+f_1x_1+\cdots+f_dx_d)^{c+d+1}}\right)dx_1\ldots dx_{d-1}\\&=&\frac{c}{(c+d)(f_0-f_d)}\left(F_c(f_1,f_2,\ldots,f_d)-F_c(f_0,f_1\ldots,f_{d-1})\right).\end{eqnarray*}The last equality is enough to extend (\ref{BF}) from $d-1$ to $d.$ 
We now apply (\ref{BF}) to the computation of $\int_{E_{d+1}}\frac{\lambda_T(dx)}{(\<f,x\>)^c}$ when $|T|=k+1$ by changing $(d,c)$ to $(k,c-k-1)$ and we get 

$$\int_{E_{d+1}}\frac{\lambda_T(dx)}{\<f,x\>^c}=\frac{k!}{(c-k)(c-k+1)\ldots(c-1)}\sum_{i\in T}\frac{1}{f_i^{c-k}\prod _{j\neq i,j\in T}(f_j-f_i)}.$$ From this calculation, statement (\ref{CP}) becomes equivalent to 

\begin{equation}\label{ECP}\sum_{\emptyset\neq T\subset\{0,\ldots,d\}}\sum_{i\in T}\frac{1}{f_i^{c+1-|T|}\prod_{j\neq i,\ j\in T}(f_j-f_i)}=f_0\ldots f_d\sum_{i=0}^d\frac{1}{f_i^{c+1}\prod_{j\neq i,\ j\in T}(f_j-f_i)}.\end{equation}
Observe now   by inversion of summations on the left--hand side that (\ref{ECP}) can be rewritten as
$$\sum_{i=0}^d\frac{1}{f_i^c}\sum_{T\ni i}\prod_{j\neq i,\ j\in T}\frac{f_i}{(f_j-f_i)}=\sum_{i=0}^d\frac{1}{f_i^c}\prod_{j\neq i}\frac{f_j}{(f_j-f_i)}.$$
Now we easily prove that for all $i=0,\ldots,d$ \begin{equation}\label{EECP}\sum_{T\ni i}\prod_{j\neq i,\ j\in T}\frac{f_i}{(f_j-f_i)}=\prod_{j\neq i}\frac{f_j}{(f_j-f_i)}.\end{equation} To see this, it is enough to prove it for $i=0.$ Denoting $X_j=f_0/(f_j-f_0),$ equality (\ref{EECP}) for $i=0$ becomes 
$\sum_{T\subset \{1,\ldots,d\}}\prod_{j\in T}X_j=\prod_{j=1}^d(1+X_j)$ which is obviously true and proves (\ref{CP}). The remainder of the theorem is plain. $\square.$ 

\vspace{4mm}
\noindent The next proposition is a remark about the weights of a face and a vertex for $ \mathcal{B}_c(a_0,\ldots,a_d)$ when this distribution exists:

\vspace{4mm}
\noindent \textbf{Proposition 4.4:}  If $B\sim \mathcal{B}_c(a_0,\ldots,a_d)$ denote $a=a_0+\cdots+a_d$ and $a'=a_{k+1}+\cdots+a_d.$ Then 
$$\Pr(B_0=\ldots=B_k=0)=\frac{\Gamma(a)\Gamma(a'+c)}{\Gamma(a+c)\Gamma(a')},\ \ \  \Pr(B=e_i)=\frac{\Gamma(a)\Gamma(a_i+c)}{\Gamma(a+c)\Gamma(a_i)}.$$

\vspace{4mm}
\noindent \textbf{Proof:} By definition $T_c(B)=T_{a+c}(X)/T_a(X)$ where $X\sim \mathcal{D}(a_0,\ldots,a_d).$ We are willing to use Property 5 of Theorem 2.1 and we consider
\begin{eqnarray}\nonumber T_c(B)(f_0,\ldots,f_0,1,\ldots,1)&=&\frac{f_0^{a-a'}}{B(a_0,\ldots,a_d)}
\int_{T_d}\frac{x_0^{a_0-1}\ldots x_{d-1}^{a_{d-1}-1}(1-x_0-\cdots-x_{d-1})^{a_d-1}dx_0\ldots dx_{d-1}}{\left((f_0-1)(x_0+\cdots+x_k)+1\right)^{a+c}}\\&&\rightarrow_{f_0\rightarrow \infty}\frac{EF}{B(a_0,\ldots,a_d)}\label{MY}\end{eqnarray}
where \begin{equation}\label{EEE}E=\int_{(0,\infty)^{k+1}}\frac{u_0^{a_0-1}\ldots u_{k}^{a_{k}-1}du_0\ldots du_{k}}{(1+u_0+\cdots+u_k)^{a+c}}=B(a_0,\ldots,a_k,a'+c)\end{equation} and $$F=\int_{T_{d-k-1}}x_{k+1}^{a_{k+1}-1}\ldots x_{d-1}^{a_{d-1}-1}(1-x_{k+1}-\cdots-x_{d-1})^{a_d-1}dx_0\ldots dx_{d-1}=B(a_{k+1},\ldots,a_d).$$ Equality (\ref{MY}) is obtained by making the change of variable $u_i=f_0x_i$ for $i=0,\ldots,k$ and taking the limit when $f_0\rightarrow \infty.$ The second equality of (\ref{EEE}) can be easily proved by applying to $A=1+u_0+\cdots+u_k$ the equality 
$$\frac{1}{A^{a+c}}=\int_{0}^{\infty}e^{-sA}s^{a+c-1}\frac{ds}{\Gamma(a+c)}.\ \square$$ 
\section{Quasi Bernoulli and Dirichlet processes.} Recall that if $(\Omega,\alpha)$ is a measure space such that $\alpha(\Omega)=a$ is finite, the  Dirichlet process with parameter $\alpha$ is a random probability $P\sim \mathcal{D}(\alpha)$ on $\Omega$ such that for any partition $(A_0,\ldots,A_d)$ of $\Omega$ then 
$$(P(A_0),\ldots,P(A_d))\sim \mathcal{D}(\alpha(A_0),\ldots,\alpha(A_d)).$$
While the term 'process' is questionable since no idea of time  is involved in this concept, it is now well ingrained in the literature; the reason is that when $\Omega$ is the interval $[0,T]$ and $\alpha$ is the Lebesgue measure, the distribution function $P\{[0,t]\}$ for $0\leq t\leq T$ is $Y(t)/Y(T)$ where $Y$ is the standard Gamma L\'evy process. A striking property of  $P\sim \mathcal{D}(\alpha)$ is that it is almost surely purely atomic: if $(V_j)_{j=1}^{\infty}$ are iid random variables on $\Omega$ with distribution $Q=\alpha/a$ there exists random weights $(W_j)_{j=1}^{\infty}$ (that is $W_j\geq 0$ and $\sum_{j=1}^{\infty}W_j=1)$ such that
$\sum_{j=1}^{\infty}W_j\delta_{V_j}\sim \mathcal{D}(\alpha).$ The exact description of the distribution of $(W_j)_{j\geq 1}$ can be found in  the fundamental paper by Ferguson [5].  A large amount of literature has followed [5]. The long paper by James, Lijoi and Prunster [7]contains a wealth of information on the Dirichlet process $P\sim \mathcal{D}(\alpha)$ and on  the distributions of the fonctionals $\int_{\Omega}f(w)P(dw)$ with numerous references. 
The present section describes the analogous  random probability $P\sim \mathcal{B}_{k}(\alpha)$ which is such that 
for any partition $(A_0,\ldots,A_d)$ of $\Omega$ we have 
\begin{equation}\label{PBER}(P(A_0),\ldots,P(A_d))\sim \mathcal{B}_k(\alpha(A_0),\ldots,\alpha(A_d)).\end{equation}
The object  $\mathcal{B}_{k}(\alpha)$ is  natural  since for $k=1$ the random probability $P=\delta_{V}$ where $V\sim \alpha/a$ satisfies (\ref{PBER}). Not surprisingly, we will see that random probabilities on $\Omega$ satisfying (\ref{PBER}) are concentrated on at most $k$ points $V_1,\ldots,V_k$ where $V_i\sim \alpha/a$, although they will not be independent as they are in the limiting case of the Dirichlet process. Needless to say, the distribution of the random weights on these atoms will not be simpler than in the limiting case.

Before stating the theorem for general $k$ we sketch the construction of  $ \mathcal{B}_{2}(\alpha)$.
We first select $V_1\sim \alpha/a.$ Having done this, with probability $1/(a+1)$ we take $V_2=V_1$ and with probability $a/(a+1)$ we choose $V_2$ independently from $V_1$ with distribution $ \alpha/a.$ Finally we take $W_1$ uniform on $(0,1)$ and $W_2=1-W_1$ and we consider the random probability
$$P=W_1\delta_{V_1}+W_2\delta_{V_2}.$$ For seeing that (\ref{PBER}) is fulfilled for $k=2$ we denote $a_i=\alpha(A_i)$ for simplicity; we observe that
the probability for which $(P(A_0),\ldots,P(A_d))$ is equal to $e_i$ where $i=0,\ldots,d$ is exactly
$$\Pr(V_1,V_2\in A_i)=\frac{1}{a+1}\frac{a_i}{a}+\frac{a}{a+1}\frac{a_i^2}{a^2}=\frac{(a_i)_2}{(a)_2}.$$ On the other hand for $i\neq j$ we have $\Pr(V_1\in A_i,\ V_2\in A_j)=\frac{a_ia_j}{a^2}.$ As a consequence the conditional distribution of $(P(A_0),\ldots,P(A_d))$ knowing that $V_1\in A_i,\ V_2\in A_j$ is $W_1e_i+W_2e_j.$ These two facts show that  (\ref{PBER}) is correct for $k=2.$

\vspace{4mm}\noindent \textbf{Theorem 5.1.} Let $(\Omega,\alpha)$ be a measure space such that $\alpha(\Omega)=a$ is finite  and let $k$ be a positive integer. We select the  random variables $(M,X,W)$ as follows:
\begin{enumerate}\item $M=(M_1,\ldots,M_k)\in \mathbb{N}^k$ are such that $M_1+2M_2+\cdots+kM_k=k$, with the Ewens distribution with parameters $k$ and $a$: 
 $$\Pr((M_1,\ldots,M_k)=(m_1,\ldots,m_k))=C(m)\frac{a^{\sum_{j=1}^km_j}}{(a)_k}$$ where $$C(m)=C(m_1,\ldots,m_k)=\frac{k!}{\prod_{j=1}^kj^{m_j}m_j!}.$$We denote $S_j=M_1+\cdots+M_j$ and $B(M)=(b_t)_{t=1}^{S_k}$ where $b_t=j$ for $S_{j-1}<t\leq S_j.$
\item We select  i.i.d. random variables $X=(X_t)_{t=1}^{S_k}$ such that $X_t\sim \alpha/a$.
\item We select $W=(W_t)_{t=1}^{S_k}\sim \mathcal{D}(B(M))$.
\item When conditioned on $M$ the random variables $X$ and $W$ are independent. 

\end{enumerate}
Then $P=\sum_{t=1}^{S_k}W_t\delta_{X_t}$ satisfies (\ref{PBER}).

\vspace{4mm}\noindent \textbf{Comments.}\begin{enumerate} \item We say that $m=(m_1,\ldots,m_k)$ is the portrait of a  permutation $\pi$ of $\{1,\ldots,k\}$  if $\pi$ has $m_j$ cycles of order $j$ for $j=1,2,\ldots, k.$ Therefore $C(m)$ is the number of permutations with portrait $m.$ For the history and the properties of the Ewens distribution,  Johnson, Kotz  and Balakrishnan [8] Chapter 41 is a good reference. Note that $m=(m_1,\ldots,m_k)$ can be seen as the coding of a partition of the integer $k$. For instance if $k=13$ and if the partition is represented by the Ferrers diagram
$$\begin{array}{ccccc}\circ&\circ&\circ&\circ&\circ\\\circ&\circ&&&\\\circ&\circ&&&\\\circ&\circ&&&\\\circ&&&&\\\circ&&&&\end{array}$$
corresponding to the partition $1+1+2+2+2+5=13$ then $m=(2,3,0,0,5,\ldots)$ where the dots mean a sequence of $8$ zeros.
The $\sum_{i=1}^km_k$  which is the height of the Ferrers diagram is equal to 6 in the example. Finally the sequence $(b_t)_{1}^{\sum_{i=1}^km_k}$ mentioned in part 1) of the theorem  is also another coding of the partition and describes the lengths  of the rows of the Ferrers diagram from below. In the above example, $(b_1,b_2,b_3,b_4,b_5,b_6)=(1,1,2,2,2,5).$
\item Let us consider the theorem for $k=3.$ In this case
$$\Pr(M=(3,0,0))=\frac{a^3}{(a)_3},\ \Pr(M=(1,1,0))=\frac{3a^2}{(a)_3},\ \Pr(M=(0,0,1))=\frac{2a}{(a)_3},$$
$$ B(3,0,0)=(1,1,1),\ B(1,1,0)=(1,2),\ B(0,0,1)=(3).$$
Therefore, \begin{itemize}\item $P=W_1\delta_{X_1}+W_2\delta_{X_2}+W_3\delta_{X_3}$ with $(W_1,W_2,W_3)\sim \mathcal{B}_3(1,1,1)$ under conditioning on $M=(3,0,0);$
\item $P=W_1\delta_{X_1}+W_2\delta_{X_2}$ with $(W_1,W_2)\sim \mathcal{B}_3(1,2)$ under conditioning on $M=(1,1,0);$
\item $P=\delta_{X_1}$ under conditioning on $M=(0,0,1).$
 \end{itemize}
\item To illustrate the notation  of Theorem 5.1, let us  come back once more to the case $k=2.$ Above, we took informally first $V_1\sim \alpha/a$, then we took $V_2$ with a mixed distribution $\frac{1}{a+1}\delta_{V_1}+\frac{a}{a+1}\frac{\alpha}{a}$ and finally we took  $P=W_1\delta_{V_1}+W_2\delta_{V_2}.$ With the new notation,  $M$ takes two values 
\begin{itemize}\item$(M_1,M_2)=(0,1)$ with probability $1/(a+1)$; in this case $X_1=V_1=V_2$, $B(0,1)=(2)$ and $P=\delta_{X_1};$
\item 
$(M_1,M_2)=(2,0)$ with probability $a/(a+1)$; in this case $B(2,0)=(1,1)$, the random probability $P$ has in general two atoms $X_1$ and $X_2$ (at least when $\alpha$ has no atoms) and they are mixed by $(W_1,W_2)=(W_1,1-W_1)\sim \mathcal{D}(1,1)$, that is to say with $W_1$  uniform on $(0,1).$\end{itemize}
\item If we consider the particular case where $\Omega=(0,1)$ and where $\alpha$ is the Lebesgue measure (therefore $a=1$) the random probability $P\sim \mathcal{B}_k(\alpha)$ on $(0,1)$ will be built according to Theorem 5.1 as follows: for creating $M$ we take a permutation $\pi$ of $\{1,\ldots,k\}$ with uniform distribution; we consider  $M=(M_1,\ldots,M_k)$  which are the number of the cycles of $\pi$ of size $1,\ldots,k$, respectively;  the sequence $M$ induces a partition $B(M)$ of the integer $k$;
we take independent random variables $X_1\ldots,X_{M_1+\cdots+M_k}$ uniformly distributed on $(0,1);$ they will be the points of discontinuity of the random distribution function $F(t)=P([0,t]);$ we take finally a Dirichlet random variable $W=(W_t)_{t=1}^{M_1+\cdots+M_k}\sim \mathcal{D}(B(M))$ and $W_t$ is the amplitude of the jump of the random process $F$ in $X_t.$

\item 
We observe that the idea of the $T_c$ transform extends well to the context of random probabilities on $\Omega.$ If $f$ is a positive measurable function on $\Omega$, if $c>0$ and if $P$ is a random probability on $\Omega$ we define $$T_c(P)(f)=\mathbb{E}\left((\int_{\Omega}f(w)P(dw))^{-c}\right)\leq \infty$$ which is finite in particular if there exists $m>0$ such that  $f(w)\geq m$ for all $w\in\Omega.$ If $P=X\sim \mathcal{D}(\alpha)$ is a Dirichlet process such that $a=\alpha(\Omega),$ formula (\ref{WW}) or [4] page 35 shows that
$$T_a(X)(f)=\mathbb{E}\left((\int_{\Omega}f(w)X(dw))^{-a}\right)=\exp-\int_{\Omega} \log f(w)\alpha(dw).$$ An interesting application of Proposition 3.1 gives the following when  $P=B\sim \mathcal{B}_k(\alpha)$ is the quasi Bernoulli process of Theorem 5.1. Denoting $\sigma_j(f)=\int_{\Omega}\frac{\alpha(dw)}{f(w)^j}$ 
we have the elegant result $$T_k(B)(f)=
\mathbb{E}\left((\int_{\Omega}f(w)B(dw))^{-k}\right)=\frac{k!}{(a)_k}\sum_{(m_1,\ldots,m_k)\in \mathbb{N}^{k}\atop  m_1+2m_2+\cdots+km_k=k}
\prod_{j=1}^k\frac{\sigma_j^{m_j}}{j^{m_j}m_j!}.$$
For instance for $k=2$ 
$$T_2(B)(f)=\frac{1}{a(a+1)}\left(\int_{\Omega}\frac{\alpha(dw)}{f(w)}\right)^2+\frac{1}{a(a+1)}\int_{\Omega}\frac{\alpha(dw)}{f(w)^2}.$$ 

\item Needless to say, the formula which is the backbone of the paper, namely
$$T_k(B)(f)=\frac{T_{a+k}(X)(f)}{T_{a}(X)(f)}$$ still holds for a Dirichlet process $X\sim \mathcal{D}(\alpha)$ and a quasi Bernoulli process $B\sim \mathcal{B}_k(\alpha).$ 
As a consequence, $X\sim (1-Y)X+YB$  holds when $Y\sim \beta(k,a)$ is independent of $B$.

\end{enumerate}

\vspace{4mm}\noindent \textbf{Proof of Theorem 5.1.} To show (\ref{PBER}) we denote $a_i=\alpha(A_i).$ We compute first the distribution of $Z=(P(A_0),\ldots,P(A_d))$ by  conditioning with respect to $M,X$. Denote by
$N_{i,j}$ the number of $X_t$ such that $S_{j-1}<t\leq S_j$ and $X_t\in A_i.$ Note that
$\sum_{i=0}^dN_{i,j}=M_j.$
The distribution of the vector $N_j=(N_{i,j})_{i=0}^d$ of $\mathbb{R}^{d+1}$ is a multinomial distribution 
$$\Pr(N_j=(n_{0,j},\ldots,n_{d,j}))=\frac{M_j!}{n_{0,j}!\cdots n_{d,j}!}\frac{a_0^{n_{0,j}}\cdots a_d^{n_{d,j}}}{a^{M_j}}$$ where $n_{0,j}+\cdots+n_{d,j}=M_j.$ Furthermore, $N_1,\ldots,N_k$ are independent. Now we introduce the following quantities
$B_i=\sum_{j=1}^kjN_{i,j}.$ They satisfy $\sum_{i=0}^{d}B_i=\sum_{i=0}^{d}\sum_{j=1}^{k}jN_{i,j}=\sum_{j=1}^{k}jM_j=k.$
We now observe that  conditionally on $M$ and $X$ we have
$Z\sim \mathcal{D}(B_0,\ldots,B_d).$
To see this we use the definition of  $Z$ for writing
\begin{eqnarray*}Z&=&\left(\sum_{t}W_t\mathbf{1}_{t; X_t\in A_0},\ldots, \sum_{t}W_t\mathbf{1}_{t; X_t\in A_d}\right)\\ &=&\left(\sum_{j=1}^k\sum_{S_{j-1}<t\leq S_j}W_t\mathbf{1}_{t; X_t\in A_0},\ldots,\sum_{j=1}^k\sum_{S_{j-1}<t\leq S_j}W_t\mathbf{1}_{t; X_t\in A_d}\right).\end{eqnarray*}
A property of the Dirichlet distribution is that if $b_i=\sum_{j=1}^{k_i}a_{ij}$ with $a_{ij}\geq 0$ and $i=0,\ldots,d,$ if $$(X_{ij})_{0\leq i\leq d,\ 1\leq j\leq k_i}\sim \mathcal{D}\left((a_{ij})_{0\leq i\leq d,\ 1\leq j\leq k_i}\right)$$ and if $Y_i=\sum_{j=1}^{k_i}X_{ij}$ then $(Y_0,\ldots,Y_d)\sim\mathcal{D}(b_0,\ldots,b_d).$ A quick way to see this is to use (\ref{CLD}). 
We apply this principle to $(X_{ij})=(W_t),$ to $k_i=k$, to $a_t=a_{ij}=j$ when $\sum_{S_{j-1}<t\leq S_j}$, and to $Y_i=P(A_i).$ We obtain:
$$Z\sim \mathcal{D}\left(\sum_{j=1}^kjN_{0,j},\ldots,\sum_{j=1}^kjN_{d,j}\right)=\mathcal{D}(B_0,\ldots,B_d).$$
The last step of the proof  removes the conditioning on $X$ and $M.$ Here in 
(\ref{FF4})  below we use the notation $\sigma_j$ introduced in (\ref{SIG}): 
\begin{eqnarray}\label{FBF}
\mathbb{E}\left(\frac{1}{\<f,Z\>^k}\right)&=&\mathbb{E}\left(\mathbb{E}\left(\frac{1}{\<f,Z\>^k}|M,X\right)\right)
=\mathbb{E}\left(\mathbb{E}\left(\frac{1}{f_0^{B_0}\ldots f_d^{B_d}}|M,X\right)\right)\\\label{FF2}
&=&\mathbb{E}\left(\mathbb{E}\left(\frac{1}{\prod_{j=1}^kf_0^{jN_{0,j}}\ldots f_d^{jN_{d,j}}}|M,X\right)\right)\\\label{FF3}
&=&\mathbb{E}\left(\prod_{j=1}^k\mathbb{E}\left(\frac{1}{f_0^{jN_{0,j}}\ldots f_d^{jN_{d,j}}}|M,X\right)\right)\\\label{FF4}&=&\mathbb{E}\left(\prod_{j=1}^k\frac{1}{a^{M_j}}\sigma_j^{M_j}\right)=\frac{k!}{(a)_k}\sum_{(m_1,\ldots,m_k)\in \mathbb{N}^{k};\atop m_1+2m_2+\cdots+km_k=k}
\prod_{j=1}^k\frac{\sigma_j^{m_j}}{j^{m_j}m_j!}.
\end{eqnarray}
Line (\ref{FBF}) comes from $Z\sim \mathcal{D}(B_0,\ldots,B_d)$ when conditioned on $(M,X)$ and from (\ref{CLD}), line 
(\ref{FF2}) comes from the definition of $B_0,\ldots,B_d,$   line 
(\ref{FF3}) comes from the independence of the $N_j$'s,  line 
(\ref{FF4}) comes from the generating function of a multinomial distribution. 
The formula (\ref{QBEX}) proves that $Z\sim \mathcal{B}_k(a_0,\ldots,a_d).$ $\square$

\section{A Markov chain on the tetrahedron} This section gives an   application of Theorems 1.1 and 4.1, does not contain new results and serves as a conclusion.  Consider the homogeneous Markov chain $(X(n))_{n\geq 0}$  valued in the convex hull $E_{d+1}$ of $(e_0,\ldots,e_d)$ 
with the following transition process: Given $X(n)\in E_{d+1}$ choose randomly a point $B(n+1)\in E_{d+1}$ such that $B(n+1)\sim \mathcal{B}_k(a_0,\ldots,a_d)$  and independently  a random number $Y_{n+1}\sim \beta(k,a).$ Now  draw the segment $(X(n),B(n+1))$ and take the point $$X(n+1)=X(n)(1-Y_{n+1})+B(n)Y_{n+1}$$ on this segment. Theorem 4.1 says that the Dirichlet distribution $\mathcal{D}(a_0,\ldots,a_d)$  is a stationary distribution for the Markov chain $(X(n))_{n\geq 0}.$ Recall the following principle (see [3] Proposition 1):

\vspace{4mm}
\noindent \textbf{Theorem 6.1:} If $E$ is a metric space and if $C$ is the set of continuous maps $f:E\rightarrow E$ 
let us fix a probability $\nu(df)$ on $C$. Let  $F_1,F_2,\ldots$ be a sequence of independent random variables on $C$ with the same distribution $\nu.$ Define $W_n=F_n\circ\ldots\circ F_2\circ F_1$ 
and $Z_n=F_1\circ\ldots\circ F_{n-1}\circ F_n.$  Suppose that almost surely $Z=\lim_n Z_n(x)$ exists in $E$ and does not depend on $x\in E$. Then 
\begin{enumerate}\item The distribution $\mu$ of $Z$ is a stationary distribution of the Markov chain $(W_n(x))_{n\geq 0}$ on $E;$ \item if $X$ and $F_1$ are independent and if $X\sim F_1(X)$ then $X\sim \mu.$ \end{enumerate}

\vspace{4mm}
\noindent  Choose $E=E_{d+1}.$ Apply Theorem 6.1 to  the distribution $\nu$ of the random map $F_1$ on $E_{d+1}$  defined by $F_1(x)=(1-Y_1)x+Y_1B(1)$ where $Y_1\sim \beta(k,a)$ and $B(1)\sim \mathcal{B}_k(a_0,\ldots,a_d)$ are independent. If the $F_n$ defined by  $F_n(x)=(1-Y_n)x+Y_nB(n)$ are independent with distribution $\nu$, clearly 
$$Z_n(x)=(\prod_{j=1}^n(1-Y_j))x+\sum_{k=1}^{n}\left(\prod_{j=1}^{k-1}(1-Y_j)\right)Y_kB(k)$$ converges almost surely to the sum of the following converging series
\begin{equation}\label{ZZ} Z=\sum_{k=1}^{\infty}\left(\prod_{j=1}^{k-1}(1-Y_j)\right)Y_kB(k)\end{equation} and therefore the hypotheses of Theorem 6.1 are met. As a consequence the Dirichlet law 
$\mathcal{D}(a_0,\ldots,a_d)$ is the unique stationary distribution of the Markov chain $(X(n))_{n\geq 0}$ and is the distribution of $Z$ defined by (\ref{ZZ}). Finally recall the definition of a perpetuity [6] on an affine space $\mathcal{A}.$ Let $\nu(df)$ be a probability on the space of affine transformations $f$ mapping  $\mathcal{A}$ into itself. We say that the probability $\mu$ on $\mathcal{A}$ is a perpetuity generated by $\nu$ if $X\sim F(X)$ when $F\sim \nu$ and $X\sim \mu$ are independent. If the conditions of Theorem 6.1 are met for $\nu$, there is exactly one perpetuity generated by  $\nu.$ Theorems 1.1, 4.1 and 4.3 say that the Dirichlet distribution is a perpetuity for the random affine map $F(x)=(1-Y)x+YB$ on the affine hyperplane  $\mathcal{A}$ of $\mathbb{R}^{d+1}$ containing $e_0,\ldots,e_d$ generated by various distributions $\nu$ of $(1-Y,YB).$
Theorem 6.1 shows that a Dirichlet process is also a perpetuity generated by the distribution of $(1-Y,YB)$ where  the set of probabilities on $\Omega$ replaces the tetrahedra with $d+1$ vertices and where $Y\sim \beta(k,a)$ is independent of the quasi Bernoulli process $\mathcal{B}_k(\alpha).$

\section{References}\vspace{4mm}\noindent[1] \textsc{Gergely Ambrus, P\'eter Kevei and Viktor V\'igh } (2012) 'The diminishing segment process'. \textit{Statist. Probab. Letters} \textbf{82} 191-195.

\vspace{4mm}\noindent[2] \textsc{Octavio Arzimendi and V\'ictor P\'erez-Abreu} (2010) 'On the non-classical infinite divisibility of power semicircle distributions.' \textit{Commun. Stoch. Anal.} \textbf{4} 161-178.

\vspace{4mm}\noindent[3] \textsc{Jean-Fran\c{c}ois Chamayou  and G\'erard Letac } (1991) 'Explicit stationary sequences for compositions of random functions and products of random matrices.' \textit{J. Theoret. Probab.} \textbf{4} 3-36

\vspace{4mm}\noindent[4] \textsc{Jean-Fran\c{c}ois Chamayou  and G\'erard Letac} (1994) 'A transient random walk on stochastic matrices with Dirichlet distributions.' \textit{Ann. Probab.} \textbf{22} 424-430.

\vspace{4mm}\noindent[5] \textsc{Thomas S. Ferguson} (1973) 'A Bayesian Analysis of Some Nonparametric Problems.' \textit{Ann. Statist.} \textbf{1} 209-230.

\vspace{4mm}\noindent[6] \textsc{Charles M. Goldie  and Ross A. Maller} (2000) 'Stability of perpetuities.' \textit{Ann. Probab.} \textbf{28} 1195-1218.

\vspace{4mm}\noindent[7] \textsc{Lancelot F. James, Antonio Lijoi  and Igor Pr\"{u}nster} (2010) 'On the posterior distribution of classes of random means' \textit{Bernoulli} \textbf{16} 155-180.

\vspace{4mm}\noindent[8] \textsc{Norman L. Johnson, Samuel Kotz and Narayanaswamy Balakrishnan} (2000) \textit{'Continuous Multivariate Distributions.'} Vol. 1, second edition, Wiley, New York.

\end{document}